\documentclass{amsart}
\usepackage{amsmath}
\usepackage{amssymb}

\newtheorem{theorem}[equation]{Theorem}

\newtheorem{prop}[equation]{Proposition}
\newtheorem{lemma}[equation]{Lemma}

\newcommand{\C}{\mathbf{C}}
\newcommand{\inv}{^{-1}}
\newcommand{\bsl}{\backslash}
\newcommand{\half}{\frac{1}{2}}

\title{On Bernstein's presentation of Iwahori-Hecke algebras and representations of split reductive groups over non-Archimedean local fields}
\author{Amritanshu Prasad}
\address{The Institute of Mathematical Sciences, CIT campus, Taramani, Chennai 600 113.}
\email{amri@ismc.res.in}


\begin{document}

\begin{abstract}
  This article gives conceptual statements and proofs relating parabolic induction and Jacquet functors on split reductive groups over a non-Archimedean local field to the associated Iwahori-Hecke algebra as tensoring from and restricting to parabolic subalgebras.
  The main tool is Bernstein's presentation of the Iwahori-Hecke algebra.
\end{abstract}

\maketitle
\pagestyle{myheadings}
\markboth{Amritanshu Prasad}{\sc {IWAHORI-HECKE ALGEBRAS}}
\section{Introduction}
\label{sec:intro}

\subsection{Background}
\label{sec:background}

The Iwahori-Hecke algebra $H$ associated to a reductive group $G(F)$ and an Iwahori subgroup $I$ is the convolution algebra of compactly-supported measures on $G(F)$ which are invariant under left and right translations in $I$ (see \S\ref{sec:definition} for details).

Iwahori and Matsumoto described a presentation of this algebra in \cite{MR32:2486} for Chevalley groups.
This is in terms of a Coxeter group $W_a$ known as the affine Weyl group.

In \cite{MR56:3196} Borel described a correspondence between irreducible admissible representations of $G(F)$ with non-zero vectors invariant under $I$ and irreducible finite dimensional $H$-modules.
He constructed an exact functor from the category of finite-dimensional $H$-modules to the category of admissible representations of $G(F)$ which maps irreducible objects to irreducible objects.

Casselman, using some techniques of Jacquet, showed that under the correspondence described by Borel, the irreducible admissible representations of $G(F)$ coming from irreducible $H$-modules are precisely those that occur as subquotients (or equivalently, subrepresentations) of the unramified principal series of $G(F)$ \cite[Proposition 2.6]{MR83a:22018}. It is implicit in Proposition 2.5 of this article that the Jacquet module of such a representation with respect to a minimal parabolic subgroup corresponds to restriction to a certain commutative subalgebra of $H$ (later identified by Bernstein).

In \cite{MR737932} Luszting introduced a new presentation of $H$, which he attributed to Bernstein.
Since then, Bernstein's presentation has played a major role in the study of Iwahori-Hecke algebras.

In \cite{MR98j:22028} Reeder gave a characterisation in terms of tensor products of the $H$-modules associated to the unramified principal series.
This means that the functor of induction from a minimal parabolic subgroup to $G(F)$ has a purely algebraic interpretation in terms of $H$.
Reeder's result was generalised by Jantzen in \cite{MR97b:22020} to induction functors for all parabolic subgroups.

\subsection{Overview}
\label{sec:overview}

In this article we show how Bernstein's presentation can be used to give conceptual statements and proofs of results relating parabolic induction and Jacquet functors in terms of tensoring and restriction on Iwahori-Hecke algebras.
\S\ref{sec:hecke_algebra} sets up notation and contains the definition of the Iwahori-Hecke algebra associated to $G(F)$ and $I$. \S\ref{sec:presentations} gives a quick overview of Bernstein's description of the Iwahori-Hecke algebra. More details may be found in \cite{MR97b:22020}. \S\ref{sec:induction_jacquet} introduces the universal unramified principal series and shows that it is a free module of rank one over $H$, based on \cite{MR1604812}. The next two sections give interpretations of parabolic induction and Jacquet functors respectively in terms of tensoring and restriction functors involving Iwahori-Hecke algebras and their parabolic subalgebras.
These results appear to be simpler than those of Reeder and Jantzen mentioned above, but are in fact equivalent, once a formula for going between left and right modules is obtained. This is done in the seventh and final section of this article.

\section{The Iwahori-Hecke algebra associated to a split reductive group}
\label{sec:hecke_algebra}

\subsection{Notation}
\label{sec:notation}

Let $F$ denote a non-Archimedean local field and let $O$ denote its ring of integers.
Fix a uniformising element in the maximal ideal of $O$ and call it $\pi$.
Let $q$ be the order of the residue field $\mathbf{F}_q=O/(\pi)$ of $F$.
Let $G$ be a connected split reductive group over $F$.
Fix a maximal split torus $A$ in $G$ and a Borel subgroup $B$ containing $A$.
Let $I$ denote the subgroup which is the pre-image of $B(\mathbf{F}_q)$ under the map $G(O)\to G(\mathbf{F}_q)$ induced by the canonical map $O\to \mathbf{F}_q$.
Then $I$ is an Iwahori subgroup of $G(F)$ in the sense of \cite[\S3.7]{MR546588}.

Let $X^*(A)$ denote lattice $\mathrm{Hom}(A,\mathbf{G}_m)$ of rational characters of $A$ and
let $X_*(A)$ denote the lattice $\mathrm{Hom}(\mathbf{G}_m,A)$ of rational cocharacters.
There is a bilinear pairing $X^*(A)\times X_*(A)\to \mathbf{Z}$ defined by requiring that for each $\lambda\in X^*(A)$, $\mu\in X_*(A)$ that $\langle \lambda, \mu \rangle=n$ when the morphism $\lambda\circ\mu:\mathbf{G}_m\to \mathbf{G}_m$ is $u\mapsto u^n$.
 
Let $\Phi$ denote the root system of $G$ with respect to $T$, viewed as a finite subset of $X^*(A)\otimes \mathbf{R}$.
$B$ determines a system $\Phi^+$ of positive roots in $\Phi$ and a set of simple roots $\Delta\subset \Phi^+$.
The Weyl group $W=N_GA/A$ is a Coxeter group, which acts on $X_*(A)$.
Let $S$ denote the set of reflections about the hyperplanes $\alpha=0$ in $X_*(A)\otimes \mathbf{R}$, where $\alpha\in \Delta$.
Then $(W,S)$ is a Coxeter system.
For $w\in W$, let $l(w)$ denote the length of the shortest word in $S$ that represents $w$.
To each $\alpha\in \Phi$, there exists a well-defined element $\check{\alpha}\in X_*(A)$ called the coroot corresponding to $\alpha$ (see e.g.,\cite[\S2.2]{MR546587}).
For $\mu\in X_*(A)$, one has $s\mu=\mu-\langle \alpha, \mu\rangle \check{\alpha}$.

The group $\tilde{W}=N_{G(F)}A/A(O)$ is called the \emph{extended affine Weyl group} of $G$. 
From the definitions, it follows that $\tilde{W}=\frac{A(F)}{A(O)}\rtimes W$.
$\tilde{W}$ is almost a Coxeter group. 
It is an extension of a finitely generated abelian group by a Coxeter group known as the \emph{affine Weyl group}. When $G$ is semi-simple and its derived group is simply connected, then $\tilde{W}$ coincides with $W_a$.

Given $\mu\in X_*(A)$ its inverse is usually denoted by $-\mu$ and therefore, it is convenient to denote by $\pi^\mu$ the element $\mu(\pi)\in A(F)$.
The map $\mu \mapsto \pi^\mu$ induces an isomorphism $X_*(A)\to\nolinebreak A(F)/A(O)$ of groups.
Identify $W$ with $N_{G(O)}A/A(O)$, and pick a section $W\to N_{G(O)}A$.
Denote also by $w$ the image of $w$ in $G(O)$.
Given $\tilde{w}=\pi^\mu w\in \tilde{W}$, with $\mu \in X_*(A)$ and $w\in W$, identify it with the element $\pi^\mu w$ of $G(F)$.

\subsection{Definition of the Iwahori-Hecke algebra}
\label{sec:definition}

Given two compactly supported locally constant complex-valued measures $\omega$ and $\tau$ on $G(F)$, their convolution $\omega*\tau$ is defined by requiring that
\begin{equation*}
  \int_{G(F)}f(x)d(\omega*\tau)(x)=\int_{G(F)\times G(F)} f(st)d(\omega\otimes\tau)(s,t)\mbox{ for all } f\in C^\infty_c(G(F)),
\end{equation*}
where $C^\infty_c(G(F))$ denotes the space of all compactly supported locally constant functions on $G(F)$.
The space of compactly supported locally constant complex-valued measures forms an associative $\mathbf{C}$-algebra with convolution as product.

If $\omega$ and $\tau$ are invariant under left and right translations in $I$, then so is $\omega*\tau$. Thus the $I$-invariant measures form a subalgebra of the convolution algebra of measures described above.
The resulting algebra $H$ is known as an \emph{Iwahori-Hecke algebra}.

Let $dg$ denote the Haar measure on $G(F)$ which gives $I$ total measure one. 
For each $\tilde{w}\in \tilde{W}$ denote by $dT_{\tilde{w}}(g)$ the measure $1_{I\tilde{w}I}(g)dg$ (treating $\tilde{w}$ as an element of $G(F)$ as described in \S\ref{sec:notation}).
Then $\{T_{\tilde{w}}|\tilde{w}\in \tilde{W}\}$ is a basis for $H$.

\subsection{Action on Iwahori-invariants}
\label{sec:action}

Given a right representation $(\pi,V)$ of $G(F)$ on $V$, the space $V^I$ of vectors in $V$ which are invariant under $I$ is endowed with the structure of a right $H$-module by the action:
\begin{equation*}
  \mathbf{v}\pi(\omega)=\int_{G(F)}\mathbf{v}\pi(g)d\omega(g) \mbox{ for each }\omega\in H.
\end{equation*}

\section{Bernstein decomposition in the Iwahori-Hecke algebra}
\label{sec:presentations}

The main results in this section are attributed to Joseph~Bernstein in \cite{MR737932}.

\subsection{The finite dimensional Hecke algebra}
\label{sec:finite_dimensional_hecke_algebra}

Let $H_0$ be the sub-algebra of $H$ consisting of the measures in $H$ which are supported on $G(O)$.
The map $G(O)\to G(\mathbf{F}_q)$ allows us to push forward a measure $\mu$ on $G(O)$ to a measure $\mu_*$ on the finite group $G(\mathbf{F}_q)$.
Let $\C[B(\mathbf{F}_q)\bsl G(\mathbf{F}_q)/B(\mathbf{F}_q)]$ denote the algebra of complex-valued measures on $G(\mathbf{F}_q)$ which are invariant under left and right translations in $B(\mathbf{F}_q)$.
It is easy to prove the following
\begin{prop}
  \label{prop:finite_dimensional}
  Push forward of measures induces an isomorphism of $\C$-algebras
  \begin{equation*}
    H_0\to \mathbf{C}[B(\mathbf{F}_q)\bsl G(\mathbf{F}_q)/B(\mathbf{F}_q)].
  \end{equation*}
\end{prop}
Thus $H_0$ is often referred to as the finite dimensional Hecke algebra of $G$.
$\mathbf{C}[B(\mathbf{F}_q)\bsl G(\mathbf{F}_q)/B(\mathbf{F}_q)]$ can be thought of as a sub-algebra of the complex group-ring of $G(\mathbf{F}_q)$.
The structure of this algebra is known, due to Iwahori \cite{MR0165016} (see also \cite[Chapter 10]{MR794307}).
Thus the following  presentation for $H_0$ is obtained:\\
Then $H_0$ is generated by $\{T_s|s\in S\}$. The relations are
\begin{equation*}
  T_s^2=q+(q-1)T_s\mbox{ for each } s\in S
\end{equation*}
and for each $s\neq t$ in $S$, there is a braid relation
\begin{equation*}
  T_s*T_t*T_s\cdots = T_t*T_s*T_t\cdots
\end{equation*}
with $m_{ij}$ terms on both sides, with $m_{ij}\geq 2$, same as those occuring in the braid relations in the presentation of the Coxeter group $W$ in terms of $S$.
$\{T_w|w\in W\}$ form a basis for $H_0$.
We will often use the well-known fact that if $w=w_1w_2\ldots w_k$, and $l(w)=l(w_1)+l(w_2)+\cdots l(w_k)$, then $T_w=T_{w_1}T_{w_2}\ldots T_{w_k}$ (this follows from \cite[Chapitre IV, \S2.4, Corollaire 1]{MR0240238}).

\subsection{Dominant cocharacters}\label{ss:lengthadd}
Let $N$ denote the unipotent radical of $B$.
Call an algebraic cocharacter $\mu \in X_*(A)$ \emph{dominant}
if $\pi^{\mu}(I\cap N)\pi^{-\mu} \subset I\cap N$.
We call $\mu$ antidominant if $-\mu$ is dominant.
The group $I$ has an \emph{Iwahori decomposition}, i.e.,
\begin{equation*}
  I=(I\cap N(F))(I\cap A(F))(I\cap \overline{N}(F)),
\end{equation*}
where $\overline{N}$ is the unipotent radical of the Borel subgroup of $G$ opposite to $B$ with respect to $A$.
It follows that
if $\mu \in X_*(A)$ is antidominant then $IwI\pi^\mu I=Iw\pi^\mu I$ for each $w$ in the finite Weyl group.
It follows that for such $\mu$, $T_{w}*T_{\pi^\mu}=T_{w\pi^\mu}$.
Similarly if $\mu$ is dominant then $T_{\pi^\mu}*T_w=T_{\pi^\mu w}$.

\subsection{The toric subalgebra}
\label{sec:toric_subalgebra}

The second subalgebra that appears in the Bernstein presentation is isomorphic to the complex group-ring $\C[X_*(A)]$.

For $\mu\in X_*(A)$ dominant, define
\begin{equation*}
  \Theta_\mu = \delta^{\half}(\pi^\mu) T_{\pi^\mu},
\end{equation*}
where $\delta^{\half}$ denotes the positive square root of the modulus function of $B(F)$.
For each such $\mu$, $\Theta_\mu$ is invertible in $H$.
If $\mu$ and $\eta$ are dominant, then the Iwahori decomposition of $I$ from \S\ref{ss:lengthadd} can be use to show that $I\pi^\mu I\pi^\eta I=I\pi^{\mu+\eta}I$, so that
\begin{equation*}
  \Theta_{\mu}*\Theta_{\eta}=\Theta_{\mu+\eta}.
\end{equation*}
Given any $\mu\in X_*(A)$, $\mu$ can be written as $\mu^+-\mu^-$, where $\mu^+$ and $\mu^-$ are dominant.
Let $R=\C[X_*(A)]$.
The map $\Theta$ defined by $\Theta:\pi^\mu \mapsto \Theta_{\mu^+}*\Theta_{\mu^-}\inv$ is an isomorphism  from $R$ onto its image in $H$ (for details see \cite[\S7]{MR737932}).

\subsection{The Bernstein relations}

Every element $h\in H$ can be written in the form $h=\Theta(r)h_0$, with $r\in R$ and $h_0\in H_0$.
Moreover, the Bernstein relations:
\begin{equation}
  \label{eq:bernstein}
  T_s*\Theta_\mu=\Theta_{s\mu}*T_s+(q-1)\frac{\Theta_{\mu}-\Theta_{s\mu}}{1-\Theta_{-\check{\alpha}}}
\end{equation}
hold for all $s\in S$, $\alpha$ the simple root corresponding to $s$ and for all $\mu \in X_*(A)$.
Note that the right hand side, although written as a fraction, does lie in the image of $R$ in $H$.
In fact, since $s\mu=\mu-\langle \alpha,\mu\rangle \check{\alpha}$,
\begin{equation*}
  \frac{\Theta_{\mu}-\Theta_{s\mu}}{1-\Theta_{-\check{\alpha}}}=\Theta_\mu+\Theta_{\mu-\check{\alpha}}+\Theta_{\mu-2\check{\alpha}}+\cdots+\Theta_{s\mu+\check{\alpha}}.
\end{equation*}

\section{The universal unramified principal series}
\label{sec:induction_jacquet}
The contents of this section are merely minor variants of results in \cite{MR1604812}.
\subsection{The universal unramified character}
\label{sec:un}
\newcommand{\chiun}{\chi_{\mathrm{un}}}
Let $\chiun:A(F)/A(O)\to R^\times$ be the character defined by
\begin{equation*}
  \chiun(\pi^\mu)=1_\mu.
\end{equation*}
Given an unramified character $\chi:A(F)/A(O)\to \C^\times$, let $\C_\chi$ denote the right $R$-module whose underlying vector space is $\C$, and where $1_\mu$ acts by $\chi(\pi^\mu)$.
Then
\begin{equation*}
  \chi = \C_\chi\otimes_R\chiun.
\end{equation*}

\subsection{The universal unramified principal series}
\label{sec:universal_unramified_principal_series}

Given an unramified character $\chi:A(F)/A(O)\to \C^\times$, let $i_{B(F)}^{G(F)}\chi$ be the $H$-module
\begin{equation*}
  \{f:G(F)/I\to \C|f(ang)=\chi(a\inv)\delta^\half(a)f(g)\}
\end{equation*}
(the equation for $f$ holding for all $a\in A(F)$, $n\in N(F)$, $g\in G(F)$).
The $H$-action is given by
\begin{equation*}
  (f\cdot \omega)(x) = \int_{G(F)} f(xg\inv)d\omega(g) \mbox{ for all } f\in i_{B(F)}^{G(F)}\chi \mbox{ and } \omega\in H.
\end{equation*}
Then $i_{B(F)}^{G(F)}\chi$ is the $H$-module of $I$-invariants in an \emph{unramified principal series} for $G(F)$.
Analogously, define $i_{B(F)}^{G(F)}\chiun$ as
\begin{equation*}
  \{f:G/I\to R|f(ang)=\chiun(a\inv)\delta^\half(a)f(g)\}
\end{equation*}
(the equation for $f$ holding for all $a\in A(F)$, $n\in N(F)$, $g\in G(F)$).
The $H$-action is exactly as in the case of $i_{B(F)}^{G(F)}\chi$.
Note that $i_{B(F)}^{G(F)}\chiun$ is a left $R$-module as well as a right $H$-module, where the $R$ and $H$-actions commute.
In the rest of this article, such an object will be called an $(R,H)$-bimodule.

$H$ acts on itself on the right by convolution. Setting $1_\mu\cdot h=\Theta_\mu * h$ for each $\mu\in X_*(A)$ and $h\in H$ gives $H$ the structure of an $(R,H)$ bimodule.

It is not difficult to verify the following
\begin{prop}
  \label{prop:isomorphism}
  The map $\Phi:\C_\chi\otimes_R i_{B(F)}^{G(F)}\chiun \to i_{B(F)}^{G(F)}\chi$ defined by
  \begin{equation*}
    \Phi(z\otimes f)\mapsto z\chi(f(g))
  \end{equation*}
  is an isomorphism of $H$-modules.
\end{prop}

\subsection{A free module of rank one}
\label{sec:free}
Let $t_w$ denote the unique element of $i_{B(F)}^{G(F)}\chiun$ that takes value $1$ at $w$, and $0$ at $v\in W$ different from $w$.
The representatives of the elements of $W$ in $G(F)$, as chosen in \S\ref{sec:notation} form coset representatives for $B(F)\bsl G(F)/I$. Hence the set $\{t_w|w\in W\}$ is a basis of $i_{B(F)}^{G(F)}\chiun$, which is a free $R$-module of rank $|W|$.
\begin{theorem}
  \label{theorem:free_one}
  The map $\Psi:H\to i_{B(F)}^{G(F)}\chiun$ defined by $\Psi(\omega)= t_1\cdot \omega$ induces an isomorphism of $(R,H)$-bimodules (with the $(R,H)$-actions as in \S\ref{sec:universal_unramified_principal_series})
\end{theorem}
Here $t_1$ is the basis vector $t_w$ with $w=1$, the identity element of $W$.
\begin{proof}
  Clearly, $\Psi$ preserves the $H$-action.
  For $\Psi$ preserve the $R$ action, it is necessary and sufficient that for each $\mu \in X_*(A)$,
  \begin{equation}
    \label{eq:check_R-action}
    1_\mu \cdot t_1 = t_1\Theta_\mu.
  \end{equation}
  It is easy to see that $t_1\Theta_\mu(w)=0$ for all $w\in W$, $w\neq 1$.
  Moreover,
  \begin{eqnarray*}
    t_1\Theta_\mu(1) & = &\delta^\half(\pi^\mu)\int_{G(F)} t_1(g\inv) 1_{I\pi^\mu I}(g) dg\\
    & = & \delta^\half(\pi^\mu)\int_{IB(F)\cap I\pi^\mu I} t_1(g\inv)dg.
  \end{eqnarray*}
  \begin{lemma}
    $IB(F)\cap I\pi^\mu I=I\pi^\mu$.
  \end{lemma}
  \begin{proof}
    Using the Iwahori decomposition for $I$, and the fact that $\mu$ is dominant, we can use the Iwahori decomposition of $I$ to write
    \begin{eqnarray*}
      I\pi^u I & = & I\pi^\mu(I\cap N(F))(I\cap A(F))\pi^{-\mu}\pi^\mu(I\cap \overline{N}(F))\\
      & = & I \pi^\mu (I\cap \overline{N}(F)).
    \end{eqnarray*}
    Suppose $\overline{n}\in I\cap \overline{N}(F)$ is such that $\pi^\mu\overline{n}\in B(F)$, then $\overline{n}\in B(F)\cap I\cap \overline{N}(F)$, which is trivial.
  \end{proof}
  It follows that
  \begin{equation*}
    (t_1\cdot \Theta_\mu )(1)=\delta^\half(\pi^\mu)t_1(\pi^{-\mu})=1_\mu=1_\mu\cdot t_1(1).
  \end{equation*}
Therefore, $\Psi$ is a morphism of $(R,H)$-bimodules.

To see that $\Psi$ is an isomorphism, note that $t_1\cdot T_w=T_w$, so that $\Psi(\Theta_\mu T_w)=1_\mu t_w$.
The above, together with the Bernstein decomposition shows that $\Psi$ is injective, as well as surjective, hence an isomorphism.
\end{proof}

\section{Induction}
\label{sec:induction}

\subsection{Induction from the Borel subgroup}
\label{sec:induction_Borel}
The following interpretation of induction in terms of tensor products follows immediately from Proposition \ref{prop:isomorphism} and Theorem \ref{theorem:free_one}

\begin{theorem}
  \label{theorem:induction_Borel}
  The map $\Psi$ induces an isomorphism of $H$-modules
  \begin{equation*}
    \C_\chi\otimes_RH \to i_{B(F)}^{G(F)}\chi.
  \end{equation*}
\end{theorem}

\subsection{$H$ is free over $R$}
\label{sec:R-module}
Recall from Theorem \ref{theorem:free_one} that $i_{B(F)}^{G(F)}\chiun$ and $H$ are isomorphic as $R$-modules, via the isomorphism $\Psi$ which identifies $t_w\in i_{B(F)}^{G(F)}\chiun$ with $T_w\in H$ for each $w\in W$.
Thus $\{t_w|w\in W\}$ generates $i_{B(F)}^{G(F)}\chiun$ as an $R$-module.
For any $r\in R$, the support of $r\cdot t_w$ lies in $BwI$.
These are disjoint for distinct $w$'s.
Therefore the $t_w$'s are linearly independent over $R$.
Thus $H$ is a free $R$-module with basis $\{T_w|w\in W\}$.

\subsection{Parabolic subgroups and subalgebras}
\label{sec:parabolic}

Given any subset $\Delta_M$ of $\Delta$ let $W_M$ denote the subgroup of $W$ generated by the set $S_M$ of simple reflections corresponding to the roots in $\Delta_M$.
The group $P=BW_M B$ is a parabolic subgroup of $G$.
The centraliser $M$ of $\{a\in A|\alpha(a)=1\mbox{ for each }\alpha\in \Delta_M\}$ is a Levi component of $P$.
If $U$ is the unipotent radical of $P$, then $P=MU$.
$A$ is a maximal torus in $M$.
$I_M=M(F)\cap I$ is an Iwahori subgroup of $M(F)$.
The Weyl group of $M$ with respect to $A$ is $W_M$.
Let $H_M$ be the Iwahori-Hecke algebra of $M(F)$ with respect to $M(F)\cap I$.
Let $\Theta_M$ be the homomorphism $\C[X_*(A)]\to H_M$ obtained by replacing $G$ with $M$ in the definition of $\Theta$ from \S\ref{sec:toric_subalgebra}.
The finite-dimensional part of $H_M$ in its Bernstein decomposition is generated by $\{T^M_s|s\in \Delta_M\}$, where $dT^M_s(m)=1_{(I\cap M(F))s(I\cap M(F))}(m)dm$, where $dm$ is a Haar measure on $M(F)$ which gives $M\cap I$ total measure one.
Use this fact to identify $H_M$ with the subalgebra of $H$ generated by $R$ and $\{T_s|s\in \Delta_M\}$ by mapping $\Theta_M(1_\mu)\mapsto \Theta_\mu$ for $\mu\in X_*(A)$ and $T^M_s\mapsto T_s$ for $s\in S_M$.

\subsection{Coset representatives}
\label{sec:representatives}
Each coset in $W_M\bsl W$ has a unique element of minimal length.
Let ${}^MW$ denote the set of such elements.
If $\Delta_M$ denotes the simple roots of $M$,
\begin{equation*}
  {}^MW=\{w\in W\;|\;w\alpha>0 \mbox{ for each }\alpha\in \Delta_M\}.
\end{equation*}
Any $w\in W$ can be written as $w''w'$ with $w''\in W_M$ and $w'\in {}^MW$
satisfying $l(w)=l(w'')+l(w')$.
Therefore $T_w=T_{w''}T_{w'}$.

\subsection{$H$ is free over $H_M$}
\label{sec:HM-module}
By \S\ref{sec:representatives}, $\{T_{w'}|w'\in {}^MW\}$ generates $H$ as an $H_M$-module.
If $\omega_M\in H_M$ then $\omega_M\cdot t_w'$ is supported on $Pw'I$.
These double cosets are disjoint for distinct $w'\in {}^MW$.
Therefore the $T_{w'}$'s are linearly independent over $H_M$.
Thus $H$ is a free $H_M$-module with basis $\{T_{w'}|w'\in {}^MW\}$. 

\subsection{A basis for tensor products}
\label{sec:basis}
Let $(\sigma,V)$ be a right $H_M$-module. Then $V\otimes_{H_M}H$ is a right $H$-module, which we call \emph{the $H$-module induced from $V$}.
Let $\{\mathbf{v}_1,\ldots,\mathbf{v}_n\}$ be a basis of $V$.
It follows from \S\ref{sec:HM-module} that  
\begin{equation*}
  \{\mathbf{v}_i\otimes T_{w'}\;|\;1\leq i\leq n, w'\in {}^MW\}
\end{equation*}
is a basis of $V\otimes_{H_M}H$.

\subsection{A basis for induced representations}
\label{sec:basis2}
Let $(\sigma,V)$ be as above.
Let $i_{P(F)}^{G(F)}V$ be the space of functions $G(F)/I\to \C$ such that
\begin{equation*}
  f(mug)d\omega_M(m)= \delta_P^\half(m)f(g)\sigma(m\inv)d\omega_M(m)
\end{equation*}
for each $m\in M(F)$, $u\in U(F)$ and $g\in G(F)$.
Here $\delta_P^\half$ is the positive square root of the modulus function of $P$.
Define a right $H$-module structure $i_{P(F)}^{G(F)}\sigma$ on $i_{P(F)}^{G(F)}V$ by
\begin{equation*}
  (i_{P(F)}^{G(F)}\sigma)(\omega)(f)(x)=\int_{G(F)} f(xg\inv)d\omega(g) \mbox{ for each }\omega\in H.
\end{equation*}
For each $w'\in {}^MW$ and $i=1,\ldots,n$ there exists a unique $t_{w',i}\in i_{P(F)}^{G(F)}V$ satisfying
\begin{equation*}
  t_{w',i}(x)=\delta_{xw'}\mathbf{v}_i \mbox{ for each } x\in {}^MW.
\end{equation*}
For each $i$, $\{t_{w',i}|w'\in{}^MW\}$ is a basis of the space of functions in $i_{P(F)}^{G(F)}V$ whose image is in the one-dimensional subspace spanned by $\mathbf{v}_i$.
It follows that 
\begin{equation*}
  \{t_{w',i}|w'\in {}^MW, i=1,\ldots,n\}
\end{equation*}
is a basis of $i_{P(F)}^{G(F)}V$.

\subsection{Tensor products and parabolic induction}
\label{parabolic_induction}
\begin{theorem}
  \label{theorem:parabolic_induction}
  The $H$-modules $(i_{P(F)}^{G(F)}\sigma,i_{P(F)}^{G(F)}V)$ and $V\otimes_{H_M}H$ are isomorphic.
\end{theorem}
\begin{proof}
  The proof is analogous to that of \cite[Proposition 2.1.2]{MR97b:22020}:
  there exists a one-dimensional right $R$-module $(\chi,\C)$ such that $V$ is a subspace of $\C\otimes_R H_M$ 
  (by Theorem \ref{theorem:induction_Borel}, $\C\otimes_R H$ is isomorphic to $i_{B(F)\cap M(F)}^{M(F)}\chi$).
  For each $1\leq i\leq n$ write 
  \begin{equation*}
    \mathbf{v}_i=\sum_{w''\in W_M} a_{w'',i}(1\otimes T_{w''}).
  \end{equation*}
  The image of $\mathbf{v}_i\otimes T_{w'}$ in $\C\otimes_R H$ is
  \begin{equation*}
    \sum_{w''\in W_M} a_{w'',i}(1\otimes T_{w''w'}).
  \end{equation*}
  Recall that the isomorphism $i_{P(F)}^{G(F)}i_{B\cap M}^M\chi\tilde{\to}i_{B(F)}^{G(F)}\chi$ is realised by $f\mapsto \tilde{f}$ where
  \begin{equation*}
    \tilde{f}(mg)=f(g)(m) \mbox{ for each } m\in M, \; g\in G.
  \end{equation*}
  Therefore, the image of $t_{w',i}$ in $i_{B(F)}^{G(F)}V$ is given by
  \begin{eqnarray*}
    \tilde{t}_{w',i}(x''x') & = & t_{w',i}(x')(x'')\\
    & = & \delta_{w'x'}\sum_{w''\in W_M} a_{w'',i}t_{w''}(x'')\\
    & = & \delta_{w'x'}\delta_{w''x''}\alpha_{w'',i}.
  \end{eqnarray*}
  This agrees with the expression for $\mathbf{v}_i\otimes T_{w'}$ obtained above.
\end{proof}

\section{Jacquet functors}
\label{sec:jacquet}
Let $P=MU$ be a standard parabolic subgroup, as in \S\ref{sec:parabolic}.
If $\pi_0$ is a left representation of $G$ on $V_0$ then the Jacquet module $r_{P(F)}^{G(F)}\pi_0$ is the representation of $M$ on the space
\begin{equation*}
  V_{0U}=V_0/V_0(U),
\end{equation*}
where 
\begin{equation*}
  V_0(U)=\mbox{linear span of}\;\{\mathbf{v}-\pi_0(n)\mathbf{v}|n\in U, \mathbf{v}\in V_0\}
\end{equation*}
given by
\begin{equation*}
r_{P(F)}^{G(F)}\pi_0(m)\mathbf{v}_U= \delta_P^{-\frac{1}{2}}(m)(\pi_0(m)\mathbf{v})_U,
\end{equation*}
where for any vector $\mathbf{x}\in V_0$, $\mathbf{x}_U$ denotes its image in $V_{0 U}$ under the quotient map.

Let $(\pi,V)$ denote the invariants under $I$ of $V_0$ let $(r_{P(F)}^{G(F)}\pi,V_U)$ denote the invariants under $M(F)\cap I$ of $(r_{P(F)}^{G(F)}\pi_0,V_{0U})$.
The $H_M$-module $(r_{P(F)}^{G(F)}\pi,V_U)$ will be called the \emph{Jacquet module of $V$ with respect to $P$ of $(\pi,V)$}.
It is a well-known that the quotient map $V_0\to V_{0U}$ induces an isomorphism of vector spaces $V\to V_U$ (this follows from \cite[Lemma 4.7]{MR56:3196} by the transitivity property of Jacquet functors, i.e., the fact that $r_{B(F)}^{G(F)}=r_{B(F)\cap M(F)}^{M(F)}\circ r_{P(F)}^{G(F)}$).

\subsection{Description of the Jacquet module with respect to $B$}
\label{sec:JacquetB}
For any dominant cocharacter $\mu$ and any $\mathbf{v}\in V$,
\begin{equation*}
  \pi(\Theta_\mu)\mathbf{v} = \delta^\half(\pi^\mu)\sum_{I/(I\cap \pi^\mu I \pi^{-\mu})}\pi_0(i\pi^\mu)\mathbf{v}.
\end{equation*}
Since $\mu$ is dominant, each coset $i(I\cap \pi^\mu I \pi^{-\mu})$ in the above sum can be represented by an element $n_i$ of $I\cap N$.
Therefore,
\begin{equation*}
  \pi(\Theta_\mu)\mathbf{v}= \delta^\half(\pi^\mu)\sum_{I/(I\cap \pi^\mu I \pi^{-\mu})}\pi_0(n_i)\pi_0(\pi^\mu) \mathbf{v}.
\end{equation*}
Since $n_i\in N$, $(\pi_0(n_i)\pi_0(\pi^\mu) \mathbf{v})_N = (\pi_0(\pi^\mu) \mathbf{v})_N$. 
Assuming that (Lemma \ref{lemma:modulus})
\begin{equation*}
  [I:I\cap \pi^\mu I \pi^{-\mu}]=\delta\inv(\pi^\mu)
\end{equation*}
we have
\begin{equation*}
  (\pi(\Theta_\mu)\mathbf{v})_N = r_{B(F)}^{G(F)}\pi (1_\mu)\mathbf{v}.
\end{equation*}
This proves
\begin{theorem}
\label{theorem:Jacquet_B}
The Jacquet module $r_{B(F)}^{G(F)}\pi$ is obtained by restricting $\pi$ to $R$ via $\Theta$.
\end{theorem}
It only remains to prove
\begin{lemma}
  \label{lemma:modulus}
  $[I:I\cap \pi^\mu I \pi^{-\mu}]=\delta\inv(\pi^\mu)$.
\end{lemma}
\begin{proof}
  Since $\mu$ is dominant, $\pi^{-\mu}I\cap \overline{N}(F)\pi^\mu\subset I$.
  Therefore, $I\cap \pi^\mu I \pi^{-\mu}\subset (I\cap A(F))(I\cap \overline{N}(F))$.
  But the Iwahori decomposition for $I$, this means that
  \begin{equation*}
    (I\cap N(F))(I\cap \pi^\mu I \pi^{-\mu})=I.
  \end{equation*}
  This means that the map
  \begin{equation*}
    (I\cap N(F))/(\pi^\mu(I\cap N)\pi^{-\mu})\to I/(I\cap \pi^\mu I \pi^{-\mu})
  \end{equation*}
  is a bijection.
  The cardinality of the left hand side of the above equation can be calculated as a quotient of volumes with respect to a left Haar measure on $B(F)$, and $\delta(\pi^{-\mu})$ appears as the answer.
\end{proof}

\subsection{Description of Jacquet modules with respect to standard parabolics}
\label{sec:JacquetP}
If $\mu$ is dominant, it is also dominant for $M$.
Therefore by \S\ref{sec:JacquetB}, $(r_{P(F)}^{G(F)}\pi(\Theta_\mu)\mathbf{v}_U)_{M\cap N}=r_{B(F)\cap M(F)}^{M(F)} r_{P(F)}^{G(F)} \pi(1_\mu)\mathbf{v}_N=(\pi(\Theta_\mu)\mathbf{v})_N$.
Therefore, $(\pi(\Theta_\mu)\mathbf{v})_U=r_{P(F)}^{G(F)}\pi(\Theta_\mu)\mathbf{v}_U$.
For brevity, write $I_M$ for $I\cap M$.
If $w\in W_M$ then the map induced by inclusion $I/(I\cap wIw\inv) \to I_M/(I_M\cap wI_Mw\inv)$ is a bijection.
Therefore
\begin{eqnarray*}
  (\pi(T_w)\mathbf{v})_U&=&\sum_{I/(I\cap wIw\inv)}(\pi_0(iw)\mathbf{v})_U\\
  & = & \sum_{I_M/(I_M\cap wI_Mw\inv)}\pi_0(iw)\mathbf{v}_U\\
  & = & r_{P(F)}^{G(F)}\pi(T_w) \mathbf{v}_U.
\end{eqnarray*}
Since the elements $T_w$, $w\in W_M$ and $\Theta_\mu$, $\mu$ dominant generate $H_M$, it follows that 
\begin{equation*}
  r_{P(F)}^{G(F)}\pi(\omega_M)\mathbf{v}_U = (\pi(\omega_M) \mathbf{v})_U
\end{equation*}
for all $\omega_M\in H_M$.

\section{Right and left modules}
\label{sec:right_left}
\subsection{The opposition formula}
\label{sec:opposition}
For a measure $\omega$ on $G(F)$ define $\omega^*$ by 
\begin{equation*}
\int_{G(F)} f(g)d\omega^*(g) = \int_{G(F)} f(g\inv)d\omega(g).
\end{equation*}
The mapping $\omega\mapsto \omega^*$ is an anti-involution for the convolution product.
A right $H$-module is viewed as a left $H$-module by writing $\omega\cdot m = m \cdot \omega^*$ for each $\omega\in H$ and each element $m$ in the module.
This is called the \emph{opposite action} of $\omega$ on $m$.
The measure $\omega^*$ is called the \emph{opposite} of $\omega$.
The opposite of $T_w$ is $T_{w\inv}$ for each $w\in \tilde{W}$.
The opposite of $\Theta_\mu$ is given by the \emph{opposition formula}
\begin{equation*}
  \Theta_\mu^*=T_{w_0}\inv \Theta_{-w_0\mu}T_{w_0},
\end{equation*}
where $w_0$ denotes the unique element in $W$ with maximal length.
It suffices to check the formula for $\mu$ dominant, in which case
\begin{equation*}
  \begin{array}{rclr}
  T_{w_0}\Theta_\mu^* & = & \delta^{\frac{1}{2}}(\pi^\mu)T_{w_0}T_{\pi^{\mu}}^* & \\
&  = & \delta^{\frac{1}{2}}(\pi^\mu)T_{w_0}T_{\pi^{{-\mu}}} & \\
&  = & \delta^{\frac{1}{2}}(\pi^\mu)T_{w_0\pi^{-\mu}} & [\mbox{by }\S\ref{ss:lengthadd}]\\
&  = & \delta^{\frac{1}{2}}(\pi^\mu)T_{\pi^{-w_0\mu}w_0} & \\
&  = & \delta^{\frac{1}{2}}(\pi^\mu)T_{\pi^{-w_0\mu}}T_{w_0} & [\mbox{by }\S\ref{ss:lengthadd}]\\
&  = & \Theta_{-w_0\mu}T_{w_0}& [\mbox{since } \delta(\pi^\mu)=\delta(\pi^{-w_0\mu})]
  \end{array}
\end{equation*}
\subsection{A formula on lengths}
\label{sec:formula_lengths}
Let $w_0'$ denote the element of maximal length in ${}^MW$.
Then for each $w''\in W_M$ the following formula holds:
\begin{equation*}
  l(w_0^{\prime-1}w''w_0')=l(w'').
\end{equation*}
To see this let $w_0''$ denote the element of maximal length in $W_M$ and write $w''=w_0''w_1''$.
Then since $w_0=w_0''w_0'$, $w_0^{\prime-1}w''w_0'=w_0\inv w_1''w_0'$, so that
\begin{eqnarray*}
  l(w_0^{\prime-1}w''w_0')&=& l(w_0\inv w_1''w_0')\\
&=&l(w_0)-l(w_1''w_0')\\
&=&l(w_0)-l(w_1'')-l(w_0')\\
&=&(l(w_0')+l(w_0''))-(l(w_0'')-l(w''))-l(w_0')\\
&=&l(w'').
\end{eqnarray*}
\subsection{A parabolic opposition formula}
\label{sec:parabolic_opposition}
The realisation of $H_M$ as a subalgebra of $H$ does not correspond to the one that comes from viewing $M$ as a closed subgroup of $G$.
In particular, the restriction of opposition to $H_M$ does not correspond to $m\mapsto m\inv$ in $M$.
This is clear when $M$ is a torus (the opposite of $\Theta_\mu$ is not $\Theta_{-\mu}$).
It is therefore necessary to distinguish between opposition in $H$ and opposition in $H_M$.
In this section $\omega_{M'}^\#$ denotes the opposite of $\omega_{M'}$ in $H_{M'}$ where $M'=w_0^{\prime-1}Mw_0'$.
As before, $\omega^*$ denotes the opposite of $\omega$ in $H$.
The opposite of $\omega_{M'}\in H_{M'}$ in $H$ is given by the \emph{parabolic opposition formula}
\begin{equation*}
  \omega_{M'}^*=T_{w_0'}^{-1}({}^{w^{\prime-1}_0}\omega_{M'}^\#)T_{w_0'}.
\end{equation*}
Suppose $\omega_{M'}=\Theta_\mu$. 
The right hand side of the formula is computed using the opposition formula for $H_{M'}$:
\begin{eqnarray*}
  T_{w_0'}\inv({}^{w_0^{\prime-1}}\Theta_\mu)^\#T_{w_0'} & = & T_{w_0'}\inv({}^{w_0^{\prime-1}}(T_{w_0^{\prime-1}w_0''w_0'}\inv\Theta_{-w_0^{\prime-1}w_0''w_0'\mu}T_{w_0^{\prime-1}w_0''w_0'}))T_{w_0'}\\
  &=& T_{w_0'}\inv T_{w_0''}\inv\Theta_{-w_0\mu}T_{w_0''}T_{w_0},
\end{eqnarray*}
which coincides with the left hand side by the opposition formula (\S\ref{sec:opposition}).
Now suppose $\omega_{M'}=T_{w_0^{\prime-1}w^{\prime\prime-1}w_0'}$ for some $w''\in W_{M}$.
The parabolic opposition formula then becomes
\begin{equation*}
  T_{w_0^{\prime-1}w''w_0'}=T_{w_0'}\inv T_{w''}T_{w_0'}.
\end{equation*}
But this equality is a direct consequence of the formula on lengths (\S\ref{sec:formula_lengths}).
Since the elements considered above generate $H_{M'}$ the parabolic opposition formula holds for all $\omega_{M'}\in H_{M'}$.
\subsection{Reeder's isomorphism}

In \cite{MR98j:22028} Reeder defines an isomorphism of left $H$-modules
\begin{equation*}
\phi:H\otimes_R \mathbf{C}_{\chi^{w_0}} \tilde{\to} (i_{B(F)}^{G(F)}\chi)^I
\end{equation*}
by $h\otimes 1\mapsto t_{w_0}h^*$.
Interpreted in terms of Theorem \ref{theorem:induction_Borel}, this is an isomorphism 
$\phi:H\otimes_R \C_{\chi^{w_0}}\to\C_{\chi\inv}\otimes_R H$, of \emph{left} H-modules (the latter is a left-module by the opposite action).
This map is well defined if and only if
\begin{equation*}
\phi(h\Theta_\mu\otimes 1)=1\otimes\Theta_{-w_0\mu}T_{w_0}h^*.
\end{equation*}
But this is immediate from the opposition formula.

\subsection{Jantzen's generalisation of Reeder's isomorphism}
\label{sec:jantzen}

Suppose $(\sigma,V)$ is a right $H_M$-module.
We know that $V\otimes_{H_M}H$ is the $H$-module obtained by parabolic induction.
Think of it as a left $H$-module under the opposite action.
In \cite[\S2.1]{MR97b:22020} Jantzen essentially shows that this module is isomorphic to the left $H$-module $H\otimes_{H_{M'}}{}^{w_0'}V^{\mathrm{op}}$.
Here $({}^{w_0'}\sigma,{}^{w_0'}V)$ is the $H_{M'}$ module obtained by composing $\sigma$ with the conjugation $x\mapsto w_0'xw_0^{\prime-1}$ mapping $M'$ to $M$.
The isomorphism $H\otimes_{H_{M'}}{}^{w_0'}V^{\mathrm{op}}\to V\otimes_{H_M}H$ is given by
\begin{equation*}
  h\otimes \mathbf{v}\mapsto \mathbf{v}\otimes T_{w_0'}h^*.
\end{equation*}
That this is a well defined homomorphism follows from the parabolic opposition formula (\S\ref{sec:parabolic_opposition}).
Once that is established, it is clear that it is an isomorphism, using the ideas from \S\ref{sec:basis}.

\subsection{Jacquet modules revisited (from the right)}
\label{sec:jacquet_right}
Combining \S\ref{sec:parabolic_opposition} with \S\ref{sec:JacquetP} (and using the notation there) gives a formula for the Jacquet modules of a right $H$-module $(\pi,V)$:
\begin{equation*}
  \mathbf{v}_{U'}r_{P'}^G\pi(\omega_{M'})=(\mathbf{v}\pi(T_{w_0'}\inv({}^{w_0^{\prime-1}}\omega_{M'})T_{w_0'}))_{U'}.
\end{equation*}

\subsection*{Acknowledgements}

It is a pleasure to thank A.~Krishnamoorthy and the students of Cochin University of Science \& Technology for their outstanding hospitality during the conference.
I am grateful to Robert~Kottwitz for his wonderful lectures on Iwahori-Hecke algebras at The University of Chicago in 2001 and to him and Thomas~Haines for some subsequent discussions, where I learned most of what is in this article.

\end{document}